\documentclass[11pt]{article}

\usepackage{amssymb,amsmath,amsthm,rotating}

\newcommand{\n}{\noindent}
\newcommand{\bb}[1]{\mathbb{#1}}
\newcommand{\cl}[1]{\mathcal{#1}}

\newcommand{\ovl}{\overline}

\theoremstyle{plain}
\newtheorem{thm}{Theorem}
\newtheorem{lem}[thm]{Lemma}

\newtheorem{pro}[thm]{Proposition}

\newtheorem*{lemP}{Lemma 5$'$}

\newtheorem{cor}[thm]{Corollary}

\theoremstyle{definition}

\begin{document}

\title{Essentially Reductive Hilbert Modules}

\author{Ronald G.~Douglas\footnote{This research was begun during a visit to 
India supported in part by the 
DST-NSF S\&T Cooperation Programme. During conferences in Chennai and 
Bangalore, the author had the opportunity to speak with W.B.\ Arveson about 
his results. The author would like to acknowledge that these conversations 
prompted this work.\newline
\n {\bf 2000 Mathematics Subject Classification.} 47B99, 47L80, 47L15.
\newline
\n {\bf Key words and phrases.} Essentially reductive Hilbert modules, 
multivariate operator theory, essentially normal operators, commuting 
weighted 
shifts. }}

\date{}
\maketitle

\abstract{Consider a Hilbert space obtained as the completion of the polynomials 
${\bb 
C}[\pmb{z}]$ in $m$-variables for which the monomials are orthogonal. If the 
commuting weighted 
shifts defined by the coordinate functions are essentially normal, then the same 
is true for their restrictions to  
invariant subspaces spanned by monomials. This generalizes the result of 
Arveson \cite{Arv03} in which the Hilbert space is the $m$-shift Hardy 
space 
$H^2_m$. 
He establishes his result for the case of finite multiplicity and shows the 
self-commutators lie in the Schatten $p$-class for $p>m$. We establish our 
result at the same level of generality. We also discuss the $K$-homology 
invariant defined in these cases.}

\setcounter{section}{-1}

\section{Introduction}\label{sec0}

\indent

The study of modules that are Hilbert spaces can be viewed as one approach to 
multivariate operator theory. While the underlying algebra could be almost 
anything, it is perhaps most natural to consider the polynomial ring ${\bb 
C}[\pmb{z}]$ or an algebra of holomorphic functions. In the case of a 
function 
algebra, such modules are called Hilbert modules and their study has been 
undertaken over the last two decades (cf.\ \cite{RP}, \cite{CG}). In this 
paper, we will 
use the terminology Hilbert module to refer to any module that is a Hilbert 
space but we will 
keep track of the hypotheses being assumed about the algebra.

A Hilbert module ${\cl M}$ is said to be essentially reductive (cf.\ 
\cite{RP}) if 
the operators $\{M_\varphi\}$ in ${\cl L}({\cl M})$ defined by module 
multiplication by elements $\varphi$ in the algebra are all essentially 
normal, 
that is, the self-commutators $[M^*_\varphi, M_\varphi] = M^*_\varphi 
M_\varphi 
- 
M_\varphi M^*_\varphi$ are in ${\cl K}({\cl M})$, the ideal of compact 
operators 
in ${\cl M}$, for all $\varphi$ in the algebra $A$. (One could also refer to 
such Hilbert modules as ``essentially normal.'') In this case, there is a 
close relationship between the algebra $A$ and the $C^*$-algebra, ${\cl 
J}({\cl 
M})$, generated by the collection $\{M_\varphi\mid \varphi\in A\}$, particularly 
the 
quotient 
$C(X_{\cl M}) = {\cl J}({\cl M})/{\cl K}({\cl M})$. In fact, the spectrum 
$X_{\cl M}$ can be identified as a subset of the maximal ideal space $M_A$ of 
the algebra $A$, if $A$ is a Banach algebra (cf.\ Thm.\ 1.6, \cite{KD}). 
Moreover, the 
$C^*$-extension defined by ${\cl J}({\cl M})$ yields an element in the odd 
$K$-homology group $K_1(X_{\cl M})$ of $X_{\cl M}$ (cf.\ \cite{BDF}) which is 
an invariant 
for the Hilbert module ${\cl M}$.

In the classical case of the Hardy and Bergman modules over the disk algebra 
$A({\bb D})$, both modules are essentially reductive as are the corresponding 
Hilbert 
modules for the Hardy and Bergman modules for the odd-dimensional spheres 
$\partial {\bb B}^n$ and balls ${\bb B}^n$. Moreover, in both cases the 
spectrum 
of the quotient $C^*$-algebras is the sphere, the boundary of ${\bb B}^n$, 
and 
the $K$-homology element is a generator for the group $K_1(\partial 
{\bb B}^n)\cong {\bb Z}$. However, for the polydisk ${\bb D}^n$, $n>1$, 
neither 
the Hardy nor 
the Bergman module is essentially reductive. More generally, one obtains an 
essentially reductive Hilbert module for strongly pseudo-convex domains in 
${\bb 
C}^n$, \cite{BdM}. In a somewhat different direction, the $m$-shift 
Hardy 
space 
$H^2_m$, which is a Hilbert module over ${\bb C}[\pmb{z}]$, 
is 
essentially reductive \cite{Arv98}.

Beyond the question of which Hilbert modules are essentially reductive, one 
can 
also ask which submodules and quotient modules are essentially reductive. In 
\cite{DM93}, 
Misra and I established by direct calculation that some quotient modules of 
the 
Hardy module for the bidisk algebra are essentially reductive and some are not. 
In 
this 
case, one can show that no nonzero submodule is essentially reductive using 
the 
fact that the coordinate functions define a pair of commuting isometries, 
both 
of infinite multiplicity. The question of essential reductivity for 
submodules and quotient modules of a given Hilbert module ${\cl M}$ is more 
likely to have an interesting answer, when ${\cl M}$  itself, is essentially 
reductive. In this note, we show that for ${\cl M}$ essentially reductive, 
either both a submodule ${\cl N}$ and the corresponding quotient module ${\cl 
M}/{\cl N}$ are essentially reductive or neither is. Moreover, we extend this 
result, which concerns short exact sequences of Hilbert modules, to longer 
resolutions 
of Hilbert modules.

In \cite{Arv03} Arveson showed that submodules of $H^2_m 
\otimes {\bb 
C}^k$, $1\le k < \infty$, in a certain class are essentially reductive and 
raised a more general 
question. His question concerned all submodules generated by homogeneous 
polynomials 
in ${\bb C}[\pmb{z}] \otimes {\bb C}^k$, $1\le k<\infty$, and he established 
essential reductivity in case the submodule is generated by monomials. Further, 
Arveson has informed me that, based on a recent result of Guo 
\cite{G}, the  question is answered in the affirmative for the 
general 
case when $m = 2$.

The action of the coordinate functions on $H^2_m\otimes {\bb C}^k$ can be 
seen 
to define commuting, contractive weighted shifts of multiplicity $k$. Our 
principal result is 
to extend Arveson's theorem to the case of general commuting weighted shifts 
so 
long as they define an essentially reductive Hilbert module over ${\bb 
C}[\pmb{z}]$. Further, 
we will 
show that our 
results  extend to the $p$-summable context which is what Arveson 
actually proves. Finally, we discuss the ${\cl K}$-homology class defined by 
this Hilbert module.

Before we begin, we want to thank the referee for pointing out a gap in our 
original proof in the $p$-summable case.

\section{Resolutions and Essential Reductivity}\label{sec1}

\indent

We begin with the result relating the behavior of submodules and their 
respective quotient modules for an essentially reductive module.

\begin{thm}\label{thm1}
Let ${\cl M}$ be an essentially reductive Hilbert module over the algebra 
$A$, 
${\cl N}$ be a submodule of ${\cl M}$ and ${\cl Q} = {\cl M}/{\cl N}$, the 
corresponding quotient module. Then ${\cl N}$ is essentially reductive iff 
${\cl 
Q}$ is.
\end{thm}

\begin{proof}
The result depends on a simple matrix calculation. For $\varphi$ in $A$ we 
consider the matrix for $M_\varphi$ relative to the decomposition ${\cl 
N}\oplus 
{\cl N}^\bot$ to obtain $\left(\begin{smallmatrix} A&B\\ 
C&D\end{smallmatrix}\right)$. Since ${\cl N}$ is invariant for $M_\varphi$ we 
have $C=0$. Moreover, the action of $\varphi$ on ${\cl N}$ defines the 
operator 
$A$, while the action of $\varphi$ on ${\cl Q}$ defines an operator unitarily 
equivalent to $D$.

Then a simple calculation shows that the matrix for $[M^*_\varphi, 
M_\varphi]$ 
relative to ${\cl M} = {\cl N}\oplus {\cl N}^\bot$ is 
$\left(\begin{smallmatrix} 
[A^*,A]-BB^*&A^*B-BD^*\\ B^*A-DB^*&[D^*,D]+B^*B\end{smallmatrix}\right)$. 
Since 
${\cl M}$ is essentially reductive, $[M^*_\varphi,M_\varphi]$ is compact and 
hence so are the operators $[A^*,A] -BB^*$ and $[D^*,D]+B^*B$. If ${\cl N}$ 
is 
essentially reductive, then $[A^*,A]$ is compact and hence $BB^*$ is compact. 
This implies $B^*B$ is compact and that $[D^*,D]$ is compact. Since this is 
true 
for every $\varphi$ in $A$, we see that ${\cl Q}$ is essentially reductive. 
The 
argument that ${\cl Q}$ essentially reductive implies that ${\cl N}$ is, 
proceeds 
similarly.
\end{proof}

The  same argument shows the theorem holds if one uses $p$-reductive instead of 
essentially reductive. (See Subsection 4 for the definition of $p$-reductive.)

Now recall that by a resolution of length one of the Hilbert module ${\cl 
M}_0$ 
over $A$, 
we mean that there are Hilbert modules ${\cl M}_1$ and ${\cl M}_2$ and 
module 
maps ${\cl M}_1 \overset{X_0}{\hbox to 25pt{\rightarrowfill}} {\cl M}_0$ and 
$M_2 \overset{X_1}{\hbox to 25pt{\rightarrowfill}} {\cl M}_1$ such that range 
$(X_0) = {\cl M}_0$, kernel $(X_1) = (0)$, and range $(X_1) =$ kernel$(X_0)$ 
(cf.\ \cite{DM2}). 
If 
$X_1$ and $X^*_0$ are isometries, then such a resolution can be seen to be  
equivalent to ${\cl M}_2$ being unitarily equivalent to a submodule of ${\cl 
M}_1$ with quotient 
module 
unitarily equivalent to ${\cl M}_0$.

\begin{thm}\label{thm2}
Consider a resolution of length one of the Hilbert module ${\cl M}_0$ over 
the 
algebra $A$:
\[
0 \longleftarrow {\cl M}_0 \overset{X_0}{\hbox to 25pt{\leftarrowfill}} {\cl 
M}_1 
\overset{X_1}{\hbox to 25pt{\leftarrowfill}} {\cl M}_2 \longleftarrow 0.
\]
If ${\cl M}_1$ and ${\cl M}_2$ are essentially reductive and $X^*_0$ is an 
isometry, then ${\cl M}_0$ is essentially reductive.
\end{thm}

\begin{proof}
We work at the level of $C^*$-algebras modulo the compacts. Fix an element 
$\varphi$ in $A$ and let $M_i$ be the operator on ${\cl L}({\cl M}_i)$ 
defined 
by module multiplication by $\varphi$ for $i=0,1,2$. Moreover, let $\pi$ 
denote 
the quotient maps ${\cl L}({\cl M}_i) \longrightarrow {\cl Q}({\cl M}_i) = 
{\cl 
L}({\cl M}_i)/{\cl K}({\cl M}_i)$ and ${\cl L}({\cl M}_i, {\cl M}_j) \to {\cl 
Q}({\cl M}_i, {\cl M}_j)$, for $0\le i,j\le 2$.

Assuming ${\cl M}_1$ is essentially reductive, we have that $\pi(M_1)$ is a 
normal element of ${\cl Q}({\cl M}_1)$. If ${\cl M}_2$ is essentially 
reductive, 
then $\pi(M_2)$ is a normal element of ${\cl Q}({\cl M}_2)$. Moreover, 
$\pi(X_1)$ intertwines $\pi(M_2)$ and $\pi(M_1)$, that is, $\pi(M_2) 
\pi(X_1) = \pi(X_1) \pi(M_1)$. Since $X_1$ is one-to-one and has 
closed range, we can write $X_1 = V_1P_1$, where $V_1$ is an isometry from 
${\cl 
M}_1$ to ${\cl M}_2$ and $P_1$ is a positive invertible operator on ${\cl 
M}_1$. 
In view of the Fuglede Theorem, the intertwining relation for $\pi(X_1)$ 
yields $\pi(M_2) \pi(V_1) = \pi(V_1) \pi(M_1)$. If we decompose 
${\cl Q}({\cl M}_1)$ using the projections $p_1 = \pi(V_1) \pi(V_1)^*$ 
and $1-p_1$, we obtain a matrix for $\pi(M_1)$ of the form 
$\left(\begin{smallmatrix} a&0\\ 0&d\end{smallmatrix}\right)$ as in the 
previous 
proof.  Since $\pi(M_1)$ is normal, we see that $d$ is normal. But $\pi(X_1)$ 
sets up a unitary equivalence between $d$ and $\pi(M_0)$ and hence 
$\pi(M_0)$ is normal. Since this is true for all $\varphi$ in $A$, we see 
that ${\cl M}_0$  is essentially reductive.
\end{proof}

If we weaken the hypotheses by not requiring $X^*_0$ to be an isometry, then 
the previous proof fails. The Fuglede Theorem requires both of the operators 
intertwined to be normal and hence we can't replace the intertwining operator 
by 
its partially isometric part. Of course, it would be enough to assume that 
$\pi(X_0)^*$ is an isometry.

We can extend these results to longer resolutions if we assume that we have a 
strong resolution, that is, if the module maps are all partial isometries (cf.\ 
\cite{DM2}). 
(Actually the last module map need not be an isometry but the others do or at 
least partial isometries modulo the compacts.) 

\begin{thm}\label{thm3}
Consider a strong resolution of finite length of the Hilbert module ${\cl 
M}_0$ 
over the algebra $A$:
\[
0 \longleftarrow {\cl M}_0 \overset{X_0}{\hbox to 25pt{\leftarrowfill}} {\cl 
M}_1 \longleftarrow \cdots \overset{X_n}{\hbox to 25pt{\leftarrowfill}} {\cl 
M}_{n+1} \longleftarrow 0.
\]
If each ${\cl M}_i$, $1\le i\le n+1$, is essentially reductive, then so is 
${\cl 
M}_0$.
\end{thm}

\begin{proof}
The proof is the same as above once one observes that at each stage one not 
only 
concludes that modulo the compacts the operators defined by module 
multiplication 
are diagonal but so also are the connecting module maps.
\end{proof}

Extending these theorems to the $p$-reductive case would involve the 
consideration of the  Fuglede Theorem in that context.

\section{Commuting Weighted Shifts---Scalar Case}\label{sec2}

\indent

We now turn to the question of establishing the essential reductivity of 
submodules of modules defined by multivariate weighted shifts. For a fixed 
positive integer $m$, let ${\bb C}[\pmb{z}]$ denote the complex polynomials 
in 
$m$-variables with $\pmb{z} = (z_1,\ldots, z_m)$. Let $\pmb{\alpha}$ be the 
multi-index $\pmb{\alpha} = (\alpha_1,\alpha_2\ldots\alpha_m)$ with each
$\alpha_i$ a 
non-negative integer, $A_m$ be the set of multi-indices, $|\pmb{\alpha}_1| = 
\alpha_1 +\ldots+ \alpha_m$, $e_i$ the multi-index with 1 in the $i$th 
position 
and zero for all other entries, and $Z^{\pmb{\alpha}}$ the monomial 
$z^{\alpha_1}_1 z^{\alpha_2}_2 \ldots z^{\alpha_m}_m$ in ${\bb C}[\pmb{z}]$. 
Let 
$\Lambda= 
\{\lambda_{\pmb{\alpha}}\}_{\pmb{\alpha}\in A_m}$ be a set of weights, $0< 
\lambda_{\pmb{\alpha}} < \infty$, and ${\cl M}_\Lambda$ be the Hilbert space 
spanned by the orthogonal set $\{\pmb{Z}^{\pmb{\alpha}}\}_{\pmb{\alpha} \in 
A_m}$ 
with $\|\pmb{Z}^{\pmb{\alpha}}\|_{{\cl M}_\Lambda} = \lambda_{\pmb{\alpha}}$. 
(This is the standard setup to define commuting weighted shifts in which the 
monomials are orthogonal.)

First, we make two basic assumptions about the set of weights. First, we assume 
that:
\begin{itemize}
\item[$(*)$] $\lambda_{\pmb{\alpha}} \ge \lambda_{\pmb{\alpha} +e_i}$ for 
$\pmb{\alpha}\in A$ and $1\le i\le m$;
\end{itemize}
which ensures that each operator $Z_i$ defined by module multiplication by 
$z_i$ 
is a contraction, $1\le i \le m$. Therefore, ${\cl M}_\Lambda$ is a Hilbert 
module over ${\bb C}[\pmb{z}]$. (Actually, in almost all of what follows, the 
assumption that the $Z_i$ are bounded is sufficient.) Moreover, one can show 
that ${\cl M}_\Lambda$ is a module for the 
algebra of functions holomorphic on any fixed polydisk of radius greater than 
one, 
although it is, in general, not a module over 
the polydisk algebra $A({\bb D}^m)$. Further, we can show that the quotient 
algebra ${\cl J}({\cl M}_\Lambda)/{\cl C}({\cl M}_\Lambda)$, where ${\cl C}({\cl 
M}_\Lambda)$ is the commutator ideal for ${\cl J}({\cl M}_\Lambda)$, is 
isometrically isomorphic to $C(X_{{\cl M}_\Lambda})$ for some compact subset 
$X_{{\cl M}_\Lambda}$ of the closed unit polydisk ${}^{cl}{\bb D}^m$ (cf.\ Thm.\ 
1.6, \cite{KD}).

Again we say that ${\cl M}_\Lambda$ is essentially reductive if the operators in 
${\cl J}({\cl M}_\Lambda)$ are
essentially normal. Our second assumption about the weight set $\Lambda$ is: 
\begin{itemize}
\item[$(**)$] $[Z_i,Z^*_j]$ is compact for all $i,j$ with $1\le i,j\le m$. 
\end{itemize} 

It is enough to assume only that $[Z_i,Z^*_i]$ is compact for $1\le i\le m$, 
since Fuglede's Theorem shows that this assumption together with the fact that 
$[Z_i,Z_j]=0$ for $1\le i,j\le m$ implies that the cross-commutators are also 
compact. We choose this form for $(**)$ to maintain parallelism with the later 
assumptions regarding $p$-summability.

In this case, since ${\cl J}({\cl M}_\Lambda)$ is irreducible, the compact 
operators ${\cl K}({\cl M}_\Lambda)$ on 
${\cl M}_\Lambda$ are contained in ${\cl J}({\cl M}_\Lambda)$, it follows that 
${\cl C}({\cl M}_\Lambda) 
= {\cl K}({\cl M}_\Lambda)$ and the quotient algebra ${\cl J}({\cl 
M}_\Lambda)/{\cl K}({\cl M}_\Lambda) \cong C(X_\Lambda)$ for some compact 
subset $X_\Lambda$ contained in ${}^{cl}{\bb D}^m$ (see Theorem \ref{thm10}).

For $B$ a subset of $A_m$, let ${\cl M}_\Lambda(B)$ be the subspace of ${\cl 
M}_\Lambda$ spanned by $\{\pmb{Z}^{\pmb{\alpha}}\}_{\pmb{\alpha}\in B}$. A 
subset $B$ of 
$A_m$ determines a submodule ${\cl M}_\Lambda(B)$ if and only if $B$ is shift 
invariant which means that $\pmb{\alpha}$ in $B$ implies $\pmb{\alpha}+e_i$ 
is 
in 
$B$ for $1\le i \le m$.  For each $i$, $1\le i \le m$, and non-negative 
integer 
$k$, set $\sum^k_i = \{\pmb{\alpha}\in A_m\mid \alpha_i = k\}$. Then ${\cl 
M}_\Lambda(\sum^k_i)$ is a reducing subspace for $Z_j$, $1\le j\le m, j\ne 
i$.

For $\pmb{\alpha}$ in $A_m$ let $B(\pmb{\alpha})$ be the subset of $A_m$ 
consisting of all $\pmb{\beta}$ satisfying $\beta_i \ge \alpha_i$, $1\le 
i\le m$. Observe that $B(\pmb{\alpha})$ is  a shift invariant subset of 
$A_m$ which is naturally isomorphic to $A_m$ and $\{\lambda_{\pmb{\beta}}\mid 
\pmb{\beta}\in 
B(\pmb{\alpha})\}$ can be identified as a weight set for $A_m$ using this 
identification.

We note that if the weight set $\Lambda$ on $A_m$ satisfies $(*)$ and $(**)$,  
then so does the weight set obtained by restricting $\Lambda$ to 
$B(\pmb{\alpha})\subset A_m$ for $\pmb{\alpha}$ in $A_m$ with $B(\pmb{\alpha})$ 
identified with $A_m$. Further, fix $i$ and $k$, $1\le i\le m$ and 
$0<k<\infty$, and identify the polynomials in ${\bb C}[\pmb{z}]$ that omit 
$z_i$ 
with the polynomials in the $(m-1)$-variables $\{z_j\}_{j\ne i}$.  Then the 
module for $A_{m-1}$ with the weight set obtained by restricting $\Lambda$ to 
$\sum^k_i$ also 
satisfies $(*)$ and $(**)$.

We note that the weight set $\lambda_{\pmb{\alpha}} = \frac{\alpha_1! 
\alpha_2! 
\ldots 
\alpha_m!}{|\pmb{\alpha}|!}$  defines the $m$-shift Hardy module $H^2_m$ and 
Arveson 
established $(**)$ in \cite{Arv98} while $(*)$ is straightforward.  

 In 
\cite{Arv03} 
Arveson showed in this case that all submodules  generated by monomials or, 
equivalently, those that 
are 
determined by a shift invariant 
subset of $A_m$ (by Proposition \ref{pro4} below), are essentially reductive. 
Our goal is 
to extend this result to 
the case of Hilbert modules defined by weighted shifts satisfying $(*)$ and 
$(**)$. Actually, Arveson establishes his result for the finite direct 
sum 
of copies of $H^2_m$ and showed that the self-commutators are in an 
appropriate Schatten $p$-class. We will do the same. 

The following result shows that the collection of submodules generated by 
shift 
invariant subsets of $A_m$, is the same as the collection of submodules 
generated by monomials.

\begin{pro}\label{pro4}
A submodule of ${\cl M}_\Lambda$ is generated by a set of monomials 
$\{Z^{\pmb{\alpha}}\}_{\pmb{\alpha}\in C}$ for $C\subset A_m$ iff it is of the 
form 
${\cl M}(B)$ for some shift invariant subset $B$ of $A_m$. Moreover, the 
generating set of monomials can be taken to be finite.
\end{pro}

\begin{proof}
If ${\cl S}$ is generated by the set 
$\{Z^{\pmb{\alpha}}\}_{\pmb{\alpha}\in C}$, then let $B$ be the shift 
invariant subset of $A_m$ generated by $C$. Then 
$\{Z^{\pmb{\alpha}}\}_{\pmb{\alpha}\in B}$ is contained in ${\cl S}$. Hence, 
${\cl M}(B)\subset {\cl S}$ and since ${\cl M}(B)$ contains 
$\{Z^{\pmb{\alpha}}\}_{\pmb{\alpha}\in C}$, we have equality. The converse 
proceeds in the same manner. Note also the proof follows from the fact that 
$B = 
\{\pmb{\alpha}\in A_m\mid Z^{\pmb{\alpha}}\in{\cl S}\}$.

The argument that the set $C$ can be taken to be finite proceeds either using 
the finite basis result for ${\bb C}[\pmb{z}]$ or the geometry of $B$.
\end{proof}

Before proceeding we need to identify a property of the weighted shifts acting 
on 
${\cl M}_\Lambda$ which follows from $(**)$. In a preliminary version of this 
paper, the conclusion of the following lemma was  assumption $(***)$ but Ken 
Davidson 
pointed out to me that it actually follows from $(**)$. We give his proof.

\begin{lem}\label{lem5}
Let $\Lambda$ be a weight  set for $A_m$, $1\le m<\infty$, satisfying $(*)$ and 
$(**)$. If $X_i$ is the edge operator   from ${\cl M}_\Lambda(\sum^k_i)$ 
to ${\cl M}_\Lambda$ defined by the action of the operator $Z_i$, then 
$X^*_iX_i$ is 
a compact operator on $M_\Lambda(\sum^k_i)$ for 
$1\le i \le m$ and $0<k<\infty$. 
\end{lem}

\begin{proof}
Fix $i$ and consider $k=0$. Let $X_{i,0}$ be the operator defined by $Z_i$ from 
${\cl M}(\sum^0_i)$ to ${\cl M}_\Lambda$. Then $X^*_{i,0}X_{i,0}$ is equal to 
the restriction of $[Z^*_i,Z_i]$ to ${\cl M}(\sum^0_i)$ since $Z^*_i|_{{\cl 
M}(\sum^0_i)} = 0$. Because $[Z^*_i,Z_i]$ is compact, then so is 
$X^*_{i,0}X_{i,0}$. Now if we consider the case of $\sum^1_i$, then the 
restriction of $[Z^*_i,Z_i]$ to ${\cl M}(\sum^1_i)$, which is compact, is the 
sum of a compact operator,  since $X^*_{i,0}X_{i,0}$ is compact, and 
$X^*_{i,1}X_{i,1}$. Thus the latter operator is compact and we can proceed 
inductively to complete the proof.
\end{proof}

Proposition \ref{pro4} shows that submodules generated by monomials have a 
geometric description, that is, are determined by shift invariant subsets of 
$A_m$. Our proofs are accomplished by decomposing the subset that determines 
the submodule into sets invariant for one 
or more of the shifts and then reducing the compactness of the self-commutators 
to 
that of the operator acting on the entire space together with the compactness of 
the 
edge operators. Similar arguments allow us to conclude that the 
cross-commutators are also compact.

\begin{thm}\label{thm5}
Let $\Lambda$ be a weight set  for $A_m$, $1\le m<\infty$,  satisfying $(*)$ and 
$(**)$. If 
$B$ is a shift invariant subset of $A_m$, then the submodule ${\cl 
M}_\Lambda(B)$ is essentially reductive.
\end{thm}

\begin{proof}
We consider first the case $m=2$  where 
the argument is more transparent since we can identify the multi-indices 
$\pmb{\alpha} = (\alpha_1,\alpha_2)$  with the integral lattice 
points in the quarter plane. We reduce  the general case to that of a 
single monomial. 
Let $\bar{\alpha}_i = 
\inf\{\alpha_i\mid (\alpha_1,\alpha_2)\in B\}$ for $i=1,2$, $\pmb{\bar 
\alpha} = 
(\bar\alpha_1,\bar\alpha_2)$, and let $\ovl{\cl M}  
= {\cl M}_\Lambda(\ovl B)$, where $\ovl B = B(\pmb{\bar \alpha})$. Since 
$\ovl 
B$ is 
shift invariant and 
contains all $(\alpha_1,\alpha_2)$ in $B$, it follows that $\ovl{\cl M}$ 
contains ${\cl M}_\Lambda(B)$. We claim that $\ovl{\cl M}/{\cl M}_\Lambda(B)$ 
is finite dimensional and hence the essential reductivity of ${\cl 
M}_\Lambda(B)$ is equivalent to that of $\ovl{\cl M}$.

In this situation, the finite dimensionality of $\ovl{\cl M}/{\cl 
M}_\Lambda(B)$ 
is equivalent to the cardinality of $\ovl{\cl B}\backslash B$ being finite. 
To see that 
the latter holds, there must exist nonnegative integers $\beta_1$ and 
$\beta_2$ such that 
$(\bar\alpha_1,\beta_2)$ and $(\beta_1,\bar\alpha_2)$ are in $B$. But then 
$\ovl 
B\backslash B$ is contained in the set $\{(\gamma_1,\gamma_2)\in A_2\mid 
\bar\alpha_1  \le \gamma_1\le \beta_1$; $\bar\alpha_2 \le \gamma_2 \le 
\beta_2\}$, which is finite.

Now we must show that the restrictions of $Z_1$ and $Z_2$ to $\ovl{\cl M}$ 
are 
essentially normal. Consider $Z_1$. Now the self-commutator of $Z_1$ on ${\cl 
M}_\Lambda$ is the direct sum of operators on the one-dimensional subspaces 
spanned by the monomials. The same is true for the restriction $Y_1$ of $Z_1$ 
to 
${\cl 
M}_\Lambda(\ovl B)$. If we set $
\ovl B = \ovl B_1 \cup \ovl B_2$, where $\ovl B_1 = \{(\gamma_1,\gamma_2)\in 
A_2 
\mid \bar\alpha_1 < \gamma_1; \bar\alpha_2 \le \gamma_2\}$ and $\ovl B_2 = 
\{(\gamma_1,\gamma_2) \in A_2\mid \bar\alpha_1 = \gamma_1; \bar\alpha_2 \le 
\gamma_2\}$, then the restrictions of the self-commutators of $Y_1$ and   
$Z_1$ to ${\cl M}_\Lambda(\ovl B_1)$ agree and hence the former is compact by 
$(**)$.

On ${\cl M}_\Lambda(\ovl B_2)$, the restrictions of the self-commutators of 
$Y_1$ and $Z_1$ agree on ${\cl M}_\Lambda$ if $\bar\alpha_1 = 0$ and hence 
again are compact by $(**)$. If $\bar\alpha_1 > 0$, then 
the restriction of the self-commutator of $Y_1$ to ${\cl M}_\Lambda(\ovl 
B_2)$ 
equals $X^*_1X_1$, where 
$X_1$ is the edge operator defined from ${\cl M}_\Lambda(\ovl B_2)$ to ${\cl 
M}_\Lambda$ by the action of $Z_1$, which is compact by Lemma 
\ref{lem5}.

Now we repeat the argument for $Z_2$ noting that the decomposition used for 
$\ovl  B$ in this case is not the same as that used above. This completes the 
proof that the self-commutators are compact for the case $m=2$.

To conclude that the cross-commutator $[Y^*_1,Y_2]$ is compact, we can either 
appeal to Fuglede's Theorem or note that the preceding analysis can be applied. 
In particular, $[Y^*_1,Y_2]$ takes a monomial $Z^{(\alpha_1,\alpha_2)}$ in the 
submodule to a multiple of $Z^{(\alpha_1-1,\alpha_2+1)}$ if the latter monomial 
is also in the submodule or to 0 otherwise. Thus on ${\cl M}_\Lambda(\ovl B_1)$, 
we obtain a 
restriction of the cross-commutator $[Z^*_1,Z_2]$ which is compact by $(**)$ or 
an edge operator on ${\cl M}_\Lambda(\ovl B_2)$ which is 
compact by Lemma \ref{lem5}. Hence, if we know that the cross-commutators on 
${\cl M}_\Lambda$ and the edge operators are compact, then the same is true for 
the restrictions $Y_1$ and $Y_2$ to the submodule.

For $m>2$ we will use induction and hence we assume the result holds for all 
$1\le m'<m$. Let $Y_m$ be the restriction of $Z_m$ to the submodule ${\cl S}$ 
generated by $\{Z^{\pmb{\alpha}^i}\}$. We show that we can reduce the 
question 
of the essential normality of $Y_m$ to the case in which all the $\alpha^i_1$ 
are 
constant. Repeating the argument, now focusing on the second component, 
allows us 
to assume not only are the $\alpha^i_1$  constant but also the $\alpha^i_2$. 
Finally, we reach the point in which all $\alpha^i_j$ are constant for 
$j=1,2,\ldots, m-1$. In this case we have to consider $Y_m$ on a submodule 
${\cl 
S}$ generated by $\{Z^{\pmb{\alpha}^i}\}$ with $\pmb{\alpha}^i = 
(\bar\alpha_1, 
\bar\alpha_2,\ldots, \bar\alpha_{m-1},\alpha^i_m)$ which equals ${\cl 
M}(B(\pmb{\alpha}))$, where $\pmb{\alpha} = (\bar\alpha_1,\ldots, 
\bar\alpha_m)$ 
with $\bar\alpha_m = \inf\{\alpha^i_m\}$. Thus we have reduced the proof to the 
case of 
a 
submodule generated by a single monomial just as for the case $m=2$. The 
argument for showing $Y_m$ is essentially normal proceeds in the same way as 
for the $m=2$ case above, by decomposing $B(\pmb{\alpha})$ into the two 
disjoint sets $\ovl B_1 = 
\{(\gamma_1,\gamma_2\ldots \gamma_n)\in A_m\mid \bar\alpha_j \le \gamma_j$, 
$1\le j< m$; $\bar\alpha_m < \gamma_m\}$ and $\ovl B_2 = 
\{(\gamma_1,\gamma_2,\ldots, \gamma_m)\in A_m\mid \bar\alpha_j \le \gamma_j$, 
$1\le j< m$; $\bar\alpha_m =\gamma_m\}$. The arguments for the compactness 
of 
the two parts of the self-commutator of $Y_m$ are also the same as those 
given above. Also, the proof of the compactness of the cross-commutator 
$[Y^*_m,Y_i]$ for $1\le i\le m-1$ proceeds in the same manner.

Now we return to the matter of reducing the general case of showing the 
essential 
normality of the restriction $Y_m$ of $Z_m$ to a submodule of  ${\cl 
S}$ 
generated by the set $\{Z^{\pmb{\alpha}^i}\}$ to the case in which the 
$\alpha^i_1$ are all constant. If $\bar\alpha_1$ is the maximum of the set 
$\alpha^i_1$, then ${\cl S}$ can be written as the orthogonal direct sum of 
the 
submodule ${\cl S}^1$ generated by $\{Z^{(\bar\alpha_1,\alpha^i_2,\ldots, 
\alpha^i_m)}\}$ and subspaces ${\cl S}^1_\gamma$, defined for $0\le 
\gamma<\bar\alpha_1$, each of which reduces $Y_m$. Thus the self-commutator 
of 
$Y_m$ is the direct sum of the self-commutator of $Z_m$ restricted to each of 
these summands. The self-commutator of the restriction of $Z_m$ to the first 
summand is the case in which all the first components are constant. To 
complete 
the reduction, we need to define the ${\cl S}^1_\gamma$ and show that the 
self-commutators of the restriction of $Y_m$ to each of them are all compact.

For each $i,1\le i\le n$, the submodule of ${\cl M}_\Lambda$ generated by 
$Z^{\pmb{\alpha}^i}$ has the form ${\cl M}(B(\pmb{\alpha}^i))$. The subspace 
${\cl S}^1_\gamma$ is spanned by the collection of monomials $N^\gamma_1 = 
\bigcup\limits^n_{i=1} \{Z^{\pmb{\beta}}\mid \pmb{\beta} \in 
B(\pmb{\alpha}^i)$, 
$\beta_1=\gamma\}$. The fact that ${\cl S}^1_\gamma$ reduces $Y_m$ follows 
from 
the fact that a monomial $Z^{\pmb{\beta}-e_m}$ is in $N^\gamma_1$, and hence 
in 
${\cl S}^1_\gamma$, if and only if it is  in ${\cl S}$. Now we can view 
$N^\gamma_1$ 
as a subset of all monomials that omit $z_1$. After identifying this set with 
the 
polynomials in $(m-1)$-variables, we obtain  the weight set by restricting 
$\Lambda$. Then noting that it satisfies $(*)$ and $(**)$, we 
can apply the 
induction hypothesis to conclude that the restriction of the self-commutator 
of 
$Y_m$ to ${\cl S}^1_\gamma$ is compact. Thus the restriction of the 
self-commutator of $Y_m$ to each ${\cl S}^1_\gamma$ is compact which 
completes 
the reduction and the proof that $Y_m$ is essentially normal. Further, the 
argument for the cross-commutators 
proceeds as above for $[Y^*_1,Y_i], 1<i\le m$, except for ${\cl S}^1_0$. Here, 
the 
argument  also
involves the fact that the cross-commutators $[Z^*_1,Z_i]$ are compact  as well 
as the compactness of the edge operator for $Z_1$ for the 
${\cl S}^1_\gamma$, $\gamma>0$, by Lemma 
\ref{lem5}.  With the next step, the reduction to the case 
in which both the first and second components, $\alpha^i_1$ and $\alpha^i_2$, 
are each all constant, we conclude that the cross-commutators $[Y^*_1,Y_i]$ and 
$[Y^*_2,Y_j]$ are compact for $1< i\le m$ and  $2< j\le m$. Hence, when we have 
completed the 
reduction, we know that all cross-commutators $[Y^*_i,Y_j]$ are compact for 
$1\le i,j\le m$.

Finally, we can repeat the argument for the restriction of each of the 
coordinate 
operators 
$\{Z_i\}_{1\le i\le m-1}$. This also enables us to conclude that all 
cross-commutators $[Y^*_i,Y_j]$ are compact. This completes the proof.
\end{proof}

\section{Commuting Weighted Shifts---Finite Multiplicity}\label{sec3}

\indent

We can extend the above result  trivially to  the case 
of 
higher 
multiplicity in one elementary situation.

\begin{cor}\label{cor6}
Let $\Lambda$ be a weight set for $A_m$, $m\ge1$, satisfying $(*)$ and $(**)$ 
and $1\le k<\infty$. If $B$ is a shift invariant subset of $A_m$, 
then 
the submodule ${\cl M}_\Lambda(B) \otimes {\bb C}^k$ of ${\cl M}_\Lambda 
\otimes 
{\bb C}^k$ is essentially reductive. Equivalently, every submodule generated 
by 
$\{Z^{\pmb{\alpha}} \otimes {\bb C}^k\mid \pmb{\alpha}\in C\}$ for some 
subset 
$C$ of $A_m$ is essentially reductive.
\end{cor}

\begin{proof}
The result follows from the theorem since ${\cl M}_\Lambda(B) \otimes {\bb 
C}^k$ 
is the direct sum of finitely many copies of ${\cl M}_\Lambda(B)$ each of 
which reduces
all of the $Z_i$.
\end{proof}

We now extend this result to general submodules of ${\cl M}_\Lambda\otimes 
{\bb 
C}^k$, $1\le k<\infty$, generated by monomials. We begin with the case $m=2$.

\begin{thm}\label{thm7}
Let $\Lambda$ be a weight set of $A_2$ satisfying $(*)$ and $(**)$  
and 
let $1\le k<\infty$. Then the submodule ${\cl S}$ generated by the set of 
monomials $\{Z^{\pmb{\alpha}^i}\otimes \pmb{x}_i\}^n_{i=1}$ for 
$\{\pmb{\alpha}^i\} \subset A_2$ and $\{\pmb{x}_i\} \subset {\bb C}^k$ is 
essentially reductive.
\end{thm}

\begin{proof}
We show that the restriction $Y_2$ of $Z_2$ to ${\cl S}$ is essentially 
normal, 
by reducing to the case in which the $\alpha^i_1$ are all equal, where 
$\pmb{\alpha}^i  = (\alpha^i_1,\alpha^i_2)$. 

This is the same argument used in the proof of Theorem \ref{thm5} with one 
additional complication which arises from the multiplicity. If $\bar\alpha_1$ 
is 
the maximum of 
the set $\{\alpha^i_1\}$, then ${\cl S}$ can be written as the orthogonal 
direct 
sum of the submodule ${\cl S}^1$ generated by 
$\{Z^{(\bar\alpha_1,\alpha^i_2)}\otimes 
\pmb{x}_i\}$ and subspaces ${\cl S}^1_\gamma$, defined for $0\le \gamma 
<\bar\alpha_1$, 
each of which reduces $Y_2$. Thus the self-commutator of $Y_2$ is the direct 
sum 
of the self-commutators of $Z_2$ restricted to each of these summands. The 
self-commutator of the restriction of $Z_2$ to the first summand is the 
reduction to the case in which all indices $(\alpha^i_1,\alpha^i_2)$ have 
constant 
first entries. To complete this reduction, we need to define the ${\cl 
S}^1_\gamma$ and show that the self-commutators of the restriction of $Z_2$ 
to 
each of them are all compact.

For each $i,1\le i\le n$, the submodule of ${\cl M}_\Lambda\otimes {\bb C}^k$ 
generated by 
$Z^{\pmb{\alpha}^i}\otimes \pmb{x}_i$ has the form ${\cl 
M}(B(\pmb{\alpha}^i)) 
\otimes (\pmb{x}_i)$, where $(\pmb{x}_i)$ denotes the subspace of ${\bb C}^k$ 
generated by the vector $\pmb{x}_i$. The subspace ${\cl S}^1_\gamma$ is 
spanned 
by the collection of monomials $N^\gamma_1 = \bigcup\limits^n_{i=1} 
\{Z^{\pmb{\beta}} \otimes \pmb{x}_i\mid \pmb{\beta} \in B(\pmb{\alpha}_i)$, 
$\beta_1=\gamma\}$. The fact that ${\cl S}^1_\gamma$ reduces $Y_2$ follows 
from 
the fact that a monomial $Z^{\pmb{\beta}-e_2}\otimes \pmb{x}_i$ is in 
$N^\gamma_1$ and 
hence in ${\cl S}^1_\gamma$ if and only if it is in ${\cl S}$. For each $0\le 
j$, let ${\cl H}^\gamma_j$ be the subspace of ${\bb C}^k$ spanned by the 
$\pmb{x}_i$ for which $Z^{(\gamma,j)} \otimes \pmb{x}_i$ is in $N^\gamma_1$. 
Then  $\{{\cl H}^\gamma_j\}$ is a strictly increasing sequence of subspaces 
of ${\bb 
C}^k$ and there exists an increasing sequence $0\le n_1 <\cdots< n_\ell 
<\infty$ such that every ${\cl H}^\gamma_j$ is equal to one of the ${\cl 
H}^\gamma_{n_i}$ and the $n_i$ are each chosen as small as possible.. Then we 
can express ${\cl S}^1_\gamma$ as the direct sum of subspaces ${\cl 
S}^1_\gamma(i)$, $1\le i\le \ell$, where ${\cl S}^1_\gamma(i)$ is the tensor 
product of the span of the monomials $\{Z^{\pmb{\beta}}\mid \beta_1 = 
\gamma$, 
$\beta_2\ge n_i\}$ with the subspace ${\cl H}^\gamma_{n_i} \cap ({\cl 
H}^\gamma_{n_{i-1}})^\bot$ of ${\bb C}^k$, where we set ${\cl H}^\gamma_{n_0} 
= (0)$. For $n_1=0$, the 
self-commutator of $Y_2$ restricted to ${\cl S}^1_\gamma(1)$ is a direct 
summand 
of the self-commutator of $Z_2\otimes I_{{\cl H}^\gamma_{n_1}}$ and hence is 
compact by $(**)$. For all $n_i>0$, the restriction of the operator $Z_2 
\otimes I_{{\cl 
H}^\gamma_{n_i}\cap ({\cl H}^\gamma_{n_{i-1}})^\bot}$ to ${\cl 
S}^1_\gamma(i)$ is compact by Lemma \ref{lem5} and hence the self-commutator is 
also 
compact. 

Now let $Y_1$ denote the restriction of $Z_1$ to ${\cl S}$ and consider the 
cross-commutator $[Y^*_1,Y_2]$ relative to the foregoing decomposition of ${\cl 
S}$. For the first term $[Y^*_1,Y_2]$ agrees with $[Z^*_1,Z_2]$ except for an 
edge operator. For all other summands, the restriction of $[Y^*_1,Y_2]$ is 
compact because both terms involve a compact edge operator. This completes the 
reduction.

Thus we have a submodule ${\cl S}$ generated by a set of monomials 
$\{Z^{\pmb{\alpha}^i} \otimes \pmb{x}_i\}$, in which $\alpha^i_1 = a_1$ 
for all 
$i$, $1\le i\le n$. The reminder of the proof is similar to what was done in the 
preceding 
paragraph. Again, ${\cl S}$ is generated by the monomials $N_1 = 
\bigcup\limits^n_{i=1} \{Z^{\pmb{\beta}} \otimes \pmb{x}_i\mid \pmb{\beta} 
\in 
B(\pmb{\alpha}^i)\}$ and for each $0\le j$, we let ${\cl H}_j$ be the subspace 
of 
${\bb C}^k$ spanned by the $\pmb{x}_i$ for which $Z^{(a_1,j)} \otimes 
\pmb{x}_i$ 
is in $N_1$. Then  $\{{\cl H}_j\}$ is an increasing sequence of subspaces of 
${\bb C}^k$ and there exists an increasing sequence $0\le n_1<n_2 <\ldots 
<n_\ell< \infty$ such that every ${\cl H}_j$ is equal to one of the ${\cl 
H}_{n_i}$ and each $n_i$ is chosen as small as possible. Then we can 
express ${\cl S}$ as the direct sum of subspaces ${\cl S}(i)$, $i\le i\le 
\ell$, 
where ${\cl S}(i)$ is the tensor product ${\cl M}(B(a_1,i)) \otimes ({\cl 
H}_{n_i} \cap ({\cl H}_{n_{i-1}})^\bot)$, again with ${\cl H}_{n_0} = (0)$. 
The 
self-commutator of the restriction of $Z_2$ to these subspaces is compact by 
Corollary \ref{cor6} since ${\cl H}_{n_i}\cap ({\cl H}_{n_{i-1}})^\bot$ is 
finite 
dimensional. The argument for cross-commutators is similar. This completes the 
proof.
\end{proof}

We now extend this result to the case $m>2$. While our argument is similar to 
that used above, it requires not only more elaborate decompositions of the 
subsets of $A_m$ but also induction on $m$.

\begin{thm}\label{thm8}
Let $\Lambda$ be a weight set for $A_m$, $1\le m<\infty$, satisfying $(*)$ 
and $(**)$ and let $1\le k<\infty$. Then the submodule ${\cl S}$ 
generated by the set of monomials $\{Z^{\pmb{\alpha}^i} \otimes 
\pmb{x}_i\}^n_{i=1}$ for $\{\pmb{\alpha}^i\} \subset A_m$ and 
$\{\pmb{x}_i\}\subset {\bb C}^k$ is essentially reductive.
\end{thm}

\begin{proof}
Fix  $m$ and assume the result holds for all $1\le m'<m$. The previous result 
fulfills the induction hypothesis.

We want to first reduce the result to the case in which the first components 
of 
the $\pmb{\alpha}^i$ are all constant. Let $\bar \alpha_1$ be the largest 
integer in the 
given set $\{\alpha^i_1\}$. First we decompose ${\cl S}$ into the orthogonal 
direct 
sum of the submodule ${\cl S}_1$ spanned by the set 
$\{Z^{(\bar\alpha_1,\alpha^i_2, 
\ldots, \alpha^i_m)} \otimes \pmb{x}_i\}$ and ${\cl S}'_1 = {\cl S}\cap ({\cl 
S}^\bot_1)$, which reduces the restrictions $Y_2,\ldots, Y_m$ of $Z_2,\ldots, 
Z_m$, respectively, to ${\cl S}$. To see this consider the collection of 
monomials 
$N_1 = \bigcup\limits^n_{i=1} \{Z^{\pmb{\beta}} \otimes \pmb{x}_i\mid 
\pmb{\beta} \in B(\pmb{\alpha}_i), \alpha^i_1 
\le\beta_1 < \alpha_1\}$. Then ${\cl S}'_1$ is 
the span of $N_1$. Now we decompose ${\cl S}'_1$ into the orthogonal direct 
sum 
of ${\cl S}'_1(\gamma)$ for $0\le \gamma < \bar\alpha_1$, where ${\cl 
S}'_1(\gamma)$ is 
the span of the monomials $\{Z^{\pmb{\beta}}\otimes \pmb{x}_i\in N_1\mid 
\beta_1 
= \gamma\}$. Each subspace ${\cl S}'_1(\gamma)$ reduces the operators, 
$Y_2,\ldots, Y_m$. Moreover, after identifying $\sum^\gamma_1$ with 
$A_{m-1}$, 
we see that ${\cl S}'_1(\gamma)$ is a submodule of ${\cl 
M}_{\Lambda^\gamma_1}\otimes {\bb C}^k$ to which the induction hypothesis 
applies, where $\Lambda^\gamma_1$ is the weight set for $A_{m-1}$ obtained 
by 
restricting $\Lambda$ to $\sum^\gamma_1$. Therefore, the self-commutators of 
the 
restrictions of $Z_2,\ldots, Z_m$ to each ${\cl S}'_1(\gamma)$ are compact. 
Hence, we can assume the first components of all the multi-indices 
$\{\pmb{\alpha}^i\}$ are the same. Again, this same decomposition can be used to 
show that the cross-commutators $[Y^*_1,Y_i]$ for $1< i\le m$ are compact.

Now starting with such a set $\{Z^{\pmb{\alpha}^i} \otimes \pmb{x}_i\}$, we 
can 
reduce the essential normality of $Z_3,\ldots, Z_m$ and the compactness of the 
cross-commutators $[Y^*_1,Y_i]$, $2\le i\le m$; $[Y^*_2,Y_j]$, $3\le j\le m$; to 
the case in which the 
 first and second components of the $\{\pmb{\alpha}^i\}$ are each 
constant. Continuing 
we eventually reduce the essential normality of the restriction of $Z_m$ to 
${\cl S}$ as well as the compactness of all cross-commutators to the case in 
which all the $\{\pmb{\alpha}^i\}$ are constant and 
then the result follows from Corollary \ref{cor6}. Thus the restriction of 
$Z_m$ 
to ${\cl S}$ is essentially normal and all cross-commutators are compact. By 
symmetry, we conclude that all the cross-commutators $[Y^*_i,Y_j]$, $1\le i,j\le 
m$, are 
compact and hence ${\cl S}$ is essentially reductive, which 
concludes 
the proof.
\end{proof}

Using the geometry of $A_m$ and the finite dimension of ${\bb C}^k$ one can show 
that every submodule ${\cl S}$ generated by a set of monomials 
$\{Z^{\pmb{\alpha}} \otimes \pmb{x}_{\pmb{\alpha}}\}_{\pmb{\alpha}\in C}$, for 
$C\subset A_m$, is finitely generated. Hence, we can 
extend Theorems \ref{thm8} and \ref{thm9} (below) to submodules generated by 
arbitrary collections of monomials.

\section{$\pmb{p}$-Summable Self Commutators}\label{sec4}

\indent

We next consider these results in the $p$-summable context. Let ${\cl C}_p$ 
denote the Schatten $p$-class (cf.\ \cite{GK}). First, we modify condition 
$(**)$  
as follows:
\medskip

\n $(**)_p$~~the cross-commutator $[Z_i,Z^*_j]$ is in ${\cl C}_p$ for $1\le 
i,j\le m$.
\medskip

In \cite{DY}, where the Berger--Shaw Theorem was generalized from single 
operators to the context of 
Hilbert modules of Krull dimension one, Hilbert modules satisfying $(**)_p$ were 
called $p$-reductive. For $(**)$ we pointed out that the Fuglede Theorem allowed 
an apparent weakening in which only the self-commutators are assumed to be 
compact. But this argument involves the Calkin algebra. There is a 
generalization to the Schatten $p$-class of the Fuglede Theorem, called the 
Fuglede--Weiss Theorem \cite{S}. However, that result doesn't seem adequate to 
reduce $(**)_p$ to assuming $p$-summability only for the self-commutators.

Note that $H^2_m$ satisfies $(**)_p$ when $p>m$. For a 
general weight set, $(**)$ often holds for smaller $p$, since it is possible for 
some of the 
$Z_i$ to be 
compact or Schatten $q$-class. (However, compare 
the remark at the end of the paragraph following the proof of Corollary 
\ref{cor12}.)
\medskip

If one examines the proofs of the previous four results, including the 
constructions in them, one sees that 
in the presence of the stronger hypothesis, one can draw  stronger 
conclusions.

\begin{thm}\label{thm9}
Let $\Lambda$ be a weight set for $A_m$, $1\le m<\infty$, satisfying $(*)$ and 
$(**)_p$  and let $1\le k<\infty$. Then for the submodule ${\cl 
S}$ 
generated by the set of monomials $\{Z^{\pmb{\alpha}} \otimes 
\pmb{x}_{\pmb{\alpha}}\}_{\pmb{\alpha}\in C}$ 
for $C\subset A_m$ and $\{\pmb{x}_{\pmb{\alpha}}\}\subset {\bb C}^k$, 
the cross-commutators of the co-ordinate functions and their adjoints lie in the 
Schatten 
$p$-class.
\end{thm}

Actually, one can show the same holds for operators defined by functions  of the 
operators which are
holomorphic on a polydisk of radius greater than one.

We omit the details of the proof of this theorem but they are precisely the same 
as 
before where the condition of 
being compact is 
replaced throughout by that of being in the Schatten $p$-class noting that Lemma 
\ref{lem5} extends to this case. This is true for both the self-commutators and 
the cross-commutators.

Actually, one can often draw sharper conclusions for the cross-commutators of 
the 
operators defined by the action of the coordinate functions and their adjoints  
on the quotient 
module ${\cl M}_\Lambda/{\cl S}$ in the presence of a slightly stronger version 
of $(**)_p$ and a strengthened Lemma \ref{lem5}, which hold  for example, for 
$H^2_m$. We will not 
attempt the maximum generality which would have to take into account 
degeneracies in the module 
action. We will use $(**)^p$  to denote the new assumption for a 
weight set $\Lambda_{\cl M}$ with the $p$ tied to $m$. Before providing the 
statements, however, we need some additional notation.

Let $\pmb{i} = \{i_1,i_2\ldots i_\ell\}$ be a subset of $\{1,2,\ldots, m\}$ so 
that $i_1 < i_2 <\ldots < i_\ell$, $[\pmb{i}] = \ell$, $I^c$ the 
complement of the $\{i_j\}$ in $\{1,\ldots, m\}$, $\pmb{k} = \{k_1,k_2\ldots 
k_\ell\}$ so that $k_j\ge 0$ for $1\le j\le \ell$, and $\sum^{\pmb{k}}_{\pmb{i}} 
=\{\pmb{\alpha}\in A_m\mid \alpha_{i_j} = k_j$, $1\le j\le \ell\}$. (Note for 
$\ell=1$, we obtain simply the $\sum^k_i$ introduced earlier.) Then 
$M_\Lambda\left(\sum^{\pmb{k}}_{\pmb{i}}\right)$ is a reducing subspace for 
$Z_p$, $p$ in $I^c$. Moreover, we can identify the polynomials in the 
$(m-\ell)$-variables, $Z_p$, for $p$ in $I^c$ with ${\bb C}[\pmb{z}]$ and obtain 
a weight set $\Lambda^{\pmb{k}}_{\pmb{i}}$ by restricting $\Lambda_{\cl M}$. Our 
strengthened assumption is:\medskip

\n $(**)^p$~~For every $\pmb{i}$ and $\pmb{k}$, the commutators of the 
restriction of $Z_\ell$ and  the compression of $Z^*_n$ to ${\cl 
M}\left(\sum^{\pmb{k}}_{\pmb{i}}\right)$ for $\ell,n$ 
in $I^c$ lies in the Schatten 
$q$-class for $q>m-[\pmb{i}]$.

\begin{lemP}
Let $\Lambda$ be a weight set for $A_m$, $1\le m <\infty$, satisfing $(*)$ and 
$(**)^p$.
For every $\pmb{i}$ and $\pmb{k}$, if  $X_p$ is the operator 
defined by 
the action of $Z_p$ 
from ${\cl M}_\Lambda\left(\sum^{\pmb{k}}_{\pmb{i}}\right)$ to ${\cl 
M}_\Lambda$, then $X^*_pX_p$ lies in the Schatten $q$-class for $q>m-[\pmb{i}]$.
\end{lemP}

Now recall that one can define the Hilbert--Samuel polynomial $p_m(z)$ (a 
polynomial in one 
variable) for a Hilbert module ${\cl M}$ over ${\bb C}[\pmb{z}]$ so long as it 
is finitely generated and the dimension of ${\cl M}/[{\bb C}_0[\pmb{z}]{\cl M}]$ 
is 
finite, where ${\bb C}_0[\pmb{z}]$ is the ideal of polynomials vanishing at 
$\pmb{0}$ 
and  [~~~]
denotes the closure in ${\cl M}$ (\cite{DY2}, cf.\ Thm.\ 4.2, \cite{Arveson2}). 
The order of 
$p_m(z)$ is said to be the dimension of ${\cl M}$. 

By analyzing the 
decompositions used in the proofs of the preceding theorems, one can show that 
the dimension of a quotient module ${\cl M}\otimes {\bb C}^k/{\cl S}$ for a 
submodule ${\cl S}$ generated by monomials is the same as the smallest 
$[\pmb{i}]$, 
where $\sum^{\pmb{k}}_{\pmb{i}}$ ranges over the blocks used in the 
decompositions in the proofs. As a 
consequence one can obtain the following result.

\begin{thm}\label{thm10}
Let $\Lambda$ be a weight set for $A_m$, $1\le m<\infty$ satisfying $(*)$ and
$(**)^p$, and let $1\le k<\infty$. If ${\cl S}$ is a submodule 
of ${\cl M}_\Lambda\otimes {\bb C}^k$ generated by a set of monomials 
$\{Z^{\pmb{a}} \otimes 
\pmb{x}_{\pmb{\alpha}}\}_{\pmb{\alpha}\in C}$ for $C\subset A_m$ and 
$\{x_{\pmb{\alpha}}\} \subset {\bb C}^k$, then the commutators of the 
operators 
defined by the coordinate functions on the quotient module ${\cl M}_\Lambda 
\otimes {\bb C}^k/{\cl S}$ and their adjoints lie in the Schatten $q$-class for 
$q>d$, where $d$ is 
the dimension of ${\cl M}_\Lambda \otimes {\bb C}_k/{\cl S}$.
\end{thm}

Since, as we pointed out above, the weight set for $H^2_m$ satisfies conditions 
$(*)$ and $(**)^p$ for $p>m$, it seems likely that the stronger 
conclusion of this theorem would hold for quotient modules $H^2_m/{\cl S}$ for 
an arbitrary submodule ${\cl S}$ generated by homogeneous polynomials, assuming 
that Arveson's conjecture is valid.

\section{$\pmb{\cl K}$-Homology Classes}\label{sec5}

\indent

Let ${\cl M}_\Lambda$ be the Hilbert module over ${\bb C}[\pmb{z}]$ defined 
by 
a weight set $\Lambda$ for $A_m$ satisfying $(*)$ and $(**)$. As we 
pointed out earlier, since 
multiplication by the coordinate functions is contractive by $(*)$, ${\cl 
M}_\Lambda$ is a bounded Hilbert  module over the algebra $A({\bb D}^m_r)$ of 
functions holomorphic on a polydisk ${\bb D}^m_r$ in ${\bb C}^m$ of radius 
$r>1$. 
Again
${\cl J}(M_\Lambda)$ denotes the $C^*$-algebra in ${\cl L}({\cl M}_\Lambda)$ 
generated by $I_{{\cl M}_\Lambda}$, $\{Z_i\}$ and the compact operators ${\cl 
K}({\cl M}_\Lambda)$, and ${\cl J}({\cl M}_\Lambda)/{\cl K}({\cl M}_\Lambda)$ 
denotes the quotient algebra. By $(**)$, the quotient algebra is commutative 
and there exists a compact metric space $X_\Lambda$ such that ${\cl J}({\cl 
M}_\Lambda)/{\cl K}({\cl M}_\Lambda) \cong C(X_\Lambda)$. Using ideas from 
the proof of Theorem 1 in \cite{KD} one can show:

\begin{thm}\label{thm11}
For $\Lambda$ a weight set for $A_m$ satisfying $(*)$ and $(**)$, 
$X_\Lambda$  can be identified as a subset of the closed polydisk ${}^{cl}{\bb 
D}^m$ so that  $Z_i$ 
corresponds to the restriction of $z_i$ to $X_\Lambda$.
\end{thm}

\begin{proof}
One obtains for every $r>1$ a homomorphism from $A({\bb D}^m_r)$ to ${\cl 
J}({\cl 
M}_\Lambda)$ and then to ${\cl J}({\cl M}_\Lambda)/{\cl K}({\cl M}_\Lambda)$. 
Thus 
we have a bounded homomorphism from $A({\bb D}^m_r)$ to $C(X_\Lambda)$. Thus 
$X_\Lambda$ can be identified as a closed subset  of ${}^{cl}{\bb D}^m_r$ for 
$r>1$. But 
the identifications for $r_1,r_2>1$ must be consistent which implies 
$X_\Lambda\subset {}^{cl}{\bb D}^m$, which completes the proof.
\end{proof}

The set $X_\Lambda$ is not an arbitrary one since it must be invariant under 
multiplication by 
$e^{i\pmb{\theta}}\equiv (e^{i\theta_1}, e^{i\theta_2}, \ldots, e^{i\theta_m})$ 
because multiplying the $m$-tuple $(Z_1,\ldots, Z_m)$ by $e^{i\pmb{\theta}}$ 
yields an $m$-tuple of operators which is unitarily equivalent to the original 
one. In the case of the 
Hardy or Bergman module over ${\bb B}^m$ or the $m$-shift Hardy module 
$H^2_m$,  $X_\Lambda$ is $\partial{\bb B}^m$. If the $(Z_1,Z_2,\ldots, Z_m)$ 
form 
a spherical contraction, (that is, $Z^*_1Z_1 +\cdots+ Z^*_mZ_m\le I$, 
then $X_\Lambda$ is contained in ${}^{cl}{\bb B}^m$. In general, it need not 
equal 
$\partial{\bb B}^m$.

\begin{cor}\label{cor12}
  A weight set $\Lambda$ for $A_m$ satisfying $(*)$ and $(**)$ determines a 
canonical element $[\Lambda]$ in ${\cl 
K}_1(X_\Lambda)$.
\end{cor}

\begin{proof}
This follows from \cite{BDF} since we have an extension
$0 \to {\cl K}({\cl M}_\Lambda) \to {\cl J}({\cl M}_\Lambda) \to C(X_\Lambda) 
\to 0$ by the Theorem.
\end{proof}

This element is not  interesting unless the ordinary homology of 
$X_\Lambda$ is nontrivial. If $X_\Lambda$ is contractible, for 
example, the closed ball ${}^{cl}{\bb D}^m$ or a point, then ${\cl 
K}_1(X_\Lambda)\cong 
(0)$ and there is no invariant. In fact, in this case, one can deform the 
$m$-tuple  of operators $\{Z_i\}$ to a commuting $m$-tuple of normal 
operators \cite{BDF}. On the other hand, if $X_\Lambda = \partial{\bb B}^m$, 
then ${\cl K}_1(\partial{\bb B}^m)\cong {\bb Z}$, and there is a non-zero 
invariant. In particular, the invariant corresponds to $-1$, giving 
$\partial{\bb B}^m$ the 
standard orientation, for $\Lambda$ the weight sets for the Bergman, Hardy or 
$m$-shift Hardy modules for ${\bb B}^m$. If we consider ${\cl 
M}_\Lambda\otimes {\bb C}^k$, then the ${\cl K}$-homology element is 
multiplied by $k$. One knows for extensions over $\partial{\bb B}^m$, that the 
${\cl 
K}$-homology element is determined by the index of the Koszul complex defined 
by a commuting $m$-tuple of generators \cite{A-S}, \cite{BDF} if such an 
$m$-tuple exists. In other 
cases, one must resort to different
measures \cite{Curto}, \cite{Arveson}. Thus, in the case of $H^2_m$, the ${\cl 
K}$-homology invariant coincides with the curvature invariant of Arveson 
\cite{Arv98}. Using the main result in \cite{DV}, one can also show if 
$[\Lambda]\ne 0$ and $\Lambda$ satisfies $(**)_p$, then $p>m$.

If ${\cl S}$ is a submodule of ${\cl M}_\Lambda \otimes {\bb C}^k$ which is 
essentially reductive, then repeating the construction in the Theorem  yields 
a closed subset $X_{\cl S}$ of $X_\Lambda$ for which ${\cl 
J}({\cl S})/{\cl K}({\cl S}) \cong C(X_{\cl S})$ and hence an element $[{\cl 
S}]$ in ${\cl K}_1(X_{\cl S})$ is defined. Similarly the quotient module ${\cl 
S}^\bot = 
{\cl M}_\Lambda \otimes {\bb C}^k/{\cl S}$ yields an element $[{\cl S}^\bot]$ 
in ${\cl  K}_1(X_{{\cl S}^\bot})$. One can show that $X_{\cl S} \cup X_{{\cl 
S}^\bot} = X_\Lambda$ and that $i^1_*([{\cl S}]) + i^2_*([{\cl S}^\bot]) = 
[\Lambda]$, where $i^1$ and $i^2$ are the inclusion maps of $X_{\cl S}$ and 
$X_{{\cl S}^\bot}$ into $X_\Lambda$, respectively. If $X_\Lambda = \partial{\bb 
B}^m$ 
and $[X_\Lambda] \ne 0$, then at least one of $[X_{\cl S}]$ and $[X_{{\cl 
S}^\bot}]$ is non-trivial and the corresponding $m$-tuple of operators defined 
by $\{Z_i\}$ 
can't be perturbed to a commuting $m$-tuple of normal operators. For $k=1$, one 
might 
conjecture that $[{\cl S}^\bot] = 0$ in this case unless ${\cl S}=(0)$.

For ${\cl S}$ a submodule of ${\cl M}_\Lambda$ generated by monomials, that 
is, the case considered in this note, one can show that $X_{{\cl S}^\bot}$ is 
the common zero set 
of the monomials generating it which is contractable. Hence, $[{\cl S}^\bot] = 
0$. This argument should 
work also for ${\cl S}$ generated by homogeneous polynomials once one knows 
that ${\cl S}$ is essentially reductive.

One can use the decompositions introduced previously in Sections \ref{sec2} and 
\ref{sec3} to draw more conclusions about 
$X_{{\cl S}^\bot}$ and $[{\cl S}^\bot]$ for ${\cl S}$ generated by monomials. 
Let ${\bb C}^m_{\deg}$ 
denote all points in ${\bb C}^m$ with at least one coordinate zero and 
$Y_\Lambda = X_\Lambda\cap {\bb C}^m_{\deg}$. Then one can show

\begin{thm}\label{thm13}
If ${\cl S}$ is a submodule of ${\cl M}_\Lambda\otimes {\bb C}_k$ generated by 
the monomials $\{Z^{\pmb{\alpha}} \otimes x_{\pmb{\alpha}}\}_{\pmb{\alpha}\in 
C}$ for $C\subset A_m$, where 
$\Lambda$ 
is a weight set for $A_m$ satisfying $(*)$ and $(**)$,  and such that 
the 
$\{\pmb{x}_{\pmb{\alpha}}\}$ span ${\bb C}^k$, then $X_{{\cl S}^\bot} \subseteq 
Y_\Lambda$. 
Hence, $i^2_*[{\cl S}^\bot] = 0$ and $i^1_*[{\cl S}] = [\Lambda]$.
\end{thm}

\begin{proof}
This argument is closely related to the one sketched for Theorem \ref{thm10}. 
One proceeds by obtaining decompositions for ${\cl S}^\bot$ analogous to those 
used in the preceding proofs for ${\cl S}$ and then noting that  the pieces 
essentially 
commute and at least one of the operators $Z_1,\ldots, Z_m$ is compact for each 
piece.
\end{proof}

If the $\{\pmb{x}_i\}$ don't span ${\bb C}^k$, then it is possible for 
$i^2_*[{\cl S}^\bot]\ne 0$. With a little more effort, one can say more. 
Again one would expect the same result to hold for quotient modules defined by 
submodules generated by 
homogeneous polynomials if Arveson's conjecture is valid.

\begin{thm}\label{thm14}
Under the same hypotheses as in Theorem \ref{thm13}, $i^2_*[{\cl S}^\bot] = 
(k-\ell)[\Lambda]$ and $i^1_*[{\cl S}] = \ell[\Lambda]$, where $\ell$ is the 
dimension of the subspace spanned by the $\{\pmb{x}_{\pmb{\alpha}}\}$ in ${\bb 
C}^k$.
\end{thm}

Note that since it is possible for $[\Lambda]=0$ in $K_1(X_\Lambda)$, these 
equations might be vacuous. However, for the weight set for $H^2_m$, we obtain 
another expression for the curvature invariant introduced by Arveson.

\end{document}